\documentclass[12pt,a4paper,draft]{article}

\usepackage[russian]{babel}
\usepackage{amsmath,amsfonts,amscd,amssymb,latexsym}
\font\teneufm=eufm10 \font\seveneufm=eufm7 \font\fiveeufm=eufm5
\newfam\frakturfam
\textfont\frakturfam=\teneufm \scriptfont\frakturfam=\seveneufm
\scriptscriptfont\frakturfam=\fiveeufm

\newtheorem{proposition}{╧Ёхфыюцхэшх}

\newtheorem{lemma}{╦хььр}
\newtheorem{theorem}{╥хюЁхьр}
\newtheorem{corollary}{╤ыхфёЄтшх}

\def\bee{\begin{eqnarray}}
\def\bes{\begin{eqnarray*}}
\def\eee{\end{eqnarray}}
\def\ees{\end{eqnarray*}}

\def\Proof{{\sl ─юърчрЄхы№ёЄтю.}\ }

\begin{document}
\noindent ╙─╩ 512.55
\pagestyle{plain}
\begin{center}{
\fontsize{12pt}{16pt}\selectfont\bf%
  \MakeUppercase{
	 ╬с ртЄюьюЁЇшчьрї ётюсюфэющ рыухсЁ√ ╦ш Ёрэур 3 эрф юсырёЄ№■ ЎхыюёЄэюёЄш}
}
\end{center}

\begin{center}

{\bf └.└. └ышьсрхт}\footnote{╩юёЄрэрщёъшщ уюёєфрЁёЄтхээ√щ яхфруюушўхёъшщ
єэштхЁёшЄхЄ шь. ╙.╤єыЄрэурчшэр,
 ╩юёЄрэрщ, 110000, ╩рчрїёЄрэ,
e-mail: {\em alialimbayev@gmail.com}},
{\bf ╨.╞. ═рєЁ√чсрхт}\footnote{┼тЁрчшщёъшщ эрЎшюэры№э√щ єэштхЁёшЄхЄ
шь. ╦.═.├єьшыхтр, └ёЄрэр, ╩рчрїёЄрэ,
e-mail: {\em nauryzbaevr@gmail.com}},
{\bf ╙.╙. ╙ьшЁсрхт}\footnote{Wayne State University,
Detroit, MI 48202, USA,
e-mail: {\em umirbaev@math.wayne.edu}

╨рсюЄр т√яюыэхэр т Ёрьърї яЁюхъЄр └╨05133009 ╚эёЄшЄєЄр ьрЄхьрЄшъш ш ьрЄхьрЄшўхёъюую ьюфхышЁютрэш  ╠╬═ ╨╩}

\end{center}

\begin{abstract}
─юърчрэю, ўЄю уЁєяяр Ёєўэ√ї ртЄюьюЁЇшчьют ётюсюфэющ рыухсЁ√ ╦ш (ш ётюсюфэющ рэЄшъюььєЄрЄштэющ рыухсЁ√) Ёрэур 3 эрф яЁюшчтюы№эющ юсырёЄ№■ ЎхыюёЄэюёЄш  ты хЄё  ётюсюфэ√ь яЁюшчтхфхэшхь уЁєяя ё юс·хфшэхээющ яюфуЁєяяющ. ╧юёЄЁюхэ яЁшьхЁ фшъюую ртЄюьюЁЇшчьр ётюсюфэющ рыухсЁ√ ╦ш (ш ётюсюфэющ рэЄшъюььєЄрЄштэющ рыухсЁ√) Ёрэур 3 эрф яЁюшчтюы№э√ь хтъышфют√ь ъюы№Ўюь. ▌ЄюЄ яЁшьхЁ  ты хЄё  рэрыюуюь шчтхёЄэюую ртЄюьюЁЇшчьр └эшър \cite{Umi07} фы  рёёюЎшрЄштэ√ї рыухсЁ.
\end{abstract}

\noindent

{\bf ╩ы■ўхт√х ёыютр:} ётюсюфэр  рыухсЁр ╦ш, ртЄюьюЁЇшчь, Ёєўэющ ртЄюьюЁЇшчь, ётюсюфэюх яЁюшчтхфхэшх уЁєяя, хтъышфютр юсырёЄ№.

\section{┬тхфхэшх}

\hspace*{\parindent}

╒юЁю°ю шчтхёЄэю \cite{Jung, Kulk}, ўЄю тёх ртЄюьюЁЇшчь√ рыухсЁ√ ьэюуюўыхэют $K[x,y]$ юЄ фтєї яхЁхьхээ√ї эрф яЁюшчтюы№э√ь яюыхь  ты ■Єё  Ёєўэ√ьш.
┴юыхх Єюую, уЁєяяр ртЄюьюЁЇшчьют рыухсЁ√ $K[x,y]$ фюяєёърхЄ ёЄЁєъЄєЁє рьры№урьшЁютрээюую ётюсюфэюую яЁюшчтхфхэш , Є.х.
$$Aut(K[x,y]) \cong A \ast_{C} B,$$
уфх $A$ -- яюфуЁєяяр рЇЇшээ√ї ртЄюьюЁЇшчьют, $B$ -- яюфуЁєяяр ЄЁхєуюы№э√ї ртЄюьюЁЇшчьют ш $C=A\cap B$.
▌ЄюЄ Ёхчєы№ЄрЄ, яю ёє∙хёЄтє, с√ы фюърчрэ т 1953 уюфє ┬. трэ фхЁ ╩рыъюь \cite{Kulk}. ┬яхЁт√х Єюўэр  ЇюЁьєышЁютър с√ыр фрэр ╚.╨. ╪рЇрЁхтшўхь \cite{Shafarevich}. ┬ ЁрсюЄх ─. ╨рщЄр \cite{Wright} фюърчрэ рэрыюу ¤Єюую Ёхчєы№ЄрЄр фы  Ёєўэ√ї ртЄюьюЁЇшчьют фтєяюЁюцфхээ√ї рыухсЁ ьэюуюўыхэют эрф яЁюшчтюы№эющ юсырёЄ№■ ЎхыюёЄэюёЄш.

└тЄюьюЁЇшчь√ ётюсюфэ√ї рёёюЎшрЄштэ√ї рыухсЁ Ёрэур 2 эрф яЁюшчтюы№э√ь яюыхь \cite{Makar-Limanov, Czerniakiewicz} ш ётюсюфэ√ї рыухсЁ ╧єрёёюэр Ёрэур 2 эрф яюыхь эєыхтющ їрЁръЄхЁшёЄшъш \cite{Poisson} Єръцх  ты ■Єё  Ёєўэ√ьш. ┴юыхх Єюую, уЁєяя√ ртЄюьюЁЇшчьют рыухсЁ√ ьэюуюўыхэют $K[x,y]$, ётюсюфэющ рёёюЎшрЄштэющ рыухсЁ√ $K\left\langle x,y\right\rangle$ ш ётюсюфэющ рыухсЁ√ ╧єрёёюэр $P\{ x,y\}$ юЄ фтєї яхЁхьхээ√ї шчюьюЁЇэ√ \cite{Makar-Limanov, Czerniakiewicz, Poisson}.
└тЄюьюЁЇшчь√ фтєяюЁюцфхээ√ї яЁртюёшььхЄЁшўэ√ї рыухсЁ эрф яЁюшчтюы№э√ь яюыхь \cite{right-symmetric} Єръцх  ты ■Єё  Ёєўэ√ьш ш уЁєяяр ртЄюьюЁЇшчьют Єръющ рыухсЁ√  фюяєёърхЄ ёЄЁєъЄєЁє рьры№урьшЁютрээюую ётюсюфэюую яЁюшчтхфхэш  \cite{NK}.

┬ 2004 уюфє ╙. ╙ьшЁсрхт ш ╚. ╪хёЄръют фюърчрыш \cite{Umi04}, ўЄю ртЄюьюЁЇшчь ═рурЄ√
$$
(x+2y(xz-y^2)+(xz-y^2)^2z,y+(xz-y^2)z,z)
$$
рыухсЁ√ ьэюуюўыхэют $K[x,y,z]$  ты хЄё  фшъшь т ёыєўрх яюы  эєыхтющ їрЁръЄхЁшёЄшъш.

╙.╙. ╙ьшЁсрхт Єръцх фюърчры \cite{Umi07}, ўЄю ртЄюьюЁЇшчь └эшър
$$
(x+z(xz-zy), y+(xz-zy)z, z)
$$
ётюсюфэющ рёёюЎшрЄштэющ рыухсЁ√ $K\left\langle x,y,z\right\rangle$ эрф яюыхь эєыхтющ їрЁръЄхЁшёЄшъш  ты хЄё  фшъшь.

┬ 1964 уюфє ╧. ╩юэ \cite{Cohn} фюърчры, ўЄю ртЄюьюЁЇшчь√ ъюэхўэю яюЁюцфхээ√ї ётюсюфэ√ї рыухсЁ ╦ш эрф яЁюшчтюы№э√ь яюыхь  ты ■Єё  Ёєўэ√ьш.
 ─ц. ╦хтшэ \cite{Lewin} юсюс∙шы ¤ЄюЄ Ёхчєы№ЄрЄ фы  °ЁрщхЁют√ї ьэюуююсЁрчшщ рыухсЁ. ═ряюьэшь, ўЄю °ЁрщхЁют√ьш  ты ■Єё  ьэюуююсЁрчш  тёхї эхрёёюЎшрЄштэ√ї рыухсЁ \cite{Kurosh}, ъюььєЄрЄштэ√ї ш рэЄшъюььєЄрЄштэ√ї рыухсЁ \cite{Shirshov54}, рыухсЁ ╦ш \cite{Shirshov53, Witt} ш ёєяхЁрыухсЁ ╦ш \cite{Mikhalev, Shtern}.
╤ыхфютрЄхы№эю, ртЄюьюЁЇшчь√ ётюсюфэ√ї эхрёёюЎшрЄштэ√ї рыухсЁ, ётюсюфэ√ї ъюььєЄрЄштэ√ї ш рэЄшъюььєЄрЄштэ√ї рыухсЁ ъюэхўэюую Ёрэур эрф яюы ьш Єръцх  ты ■Єё  Ёєўэ√ьш.

┬ ЁрсюЄрї \cite{Kukin89, Kukin91} с√ыш шёёыхфютрэ√ уЁєяя√ ртЄюьюЁЇшчьют ъюэхўэю яюЁюцфхээ√ї ётюсюфэ√ї рыухсЁ эрф ъюы№Ўюь уыртэ√ї шфхрыют. ┬ ўрёЄэюёЄш юфшэ шч Ёхчєы№ЄрЄют (╥хюЁхьр 3, \cite{Kukin91}) уырёшЄ, ўЄю ртЄюьюЁЇшчь√ ётюсюфэ√ї эхрёёюЎшрЄштэ√ї рыухсЁ эрф ъюы№Ўюь уыртэ√ї шфхрыют  ты ■Єё  Ёєўэ√ьш. ╬фэръю т ЁрсюЄх \cite{AU17} яЁштхфхэ яЁшьхЁ фшъюую ртЄюьюЁЇшчьр
$$
\eta=(x_1+x_2(zx_1-x_2^2)+(zx_1-x_2^2)x_2+z(zx_1-x_2^2)^2,\\x_2+z(zx_1-x_2^2))
$$
 ётюсюфэющ эхрёёюЎшрЄштэющ рыухсЁ√ ш ётюсюфэющ ъюььєЄрЄштэющ рыухсЁ√ Ёрэур фтр эрф хтъышфют√ь ъюы№Ўюь. ▌ЄюЄ ртЄюьюЁЇшчь  ты хЄё  юсюс∙хэшхь ртЄюьюЁЇшчьр ═рурЄ√.
╠хЄюф яюёЄЁюхэш  ¤Єюую ртЄюьюЁЇшчьр эх яЁюїюфшЄ фы  рэЄшъюььєЄрЄштэ√ї рыухсЁ, Єръ ъръ ътрфЁрЄ ¤ыхьхэЄр т рэЄшъюььєЄрЄштэ√ї рыухсЁрї Ёртхэ эєы■. ╬сюс∙хэшх ьхЄюфют яюёЄЁюхэш  ртЄюьюЁЇшчьют ═рурЄ√ \cite{Nagata} ш └эшър \cite{Cohn72} яючтюышыю т эрёЄю ∙хщ ЁрсюЄх яюёЄЁюшЄ№ фшъшщ ртЄюьюЁЇшчь
$$
\delta=(x_1+[zx_1-[x_2,x_3],x_3], x_2+z(zx_1-[x_2, x_3]), x_3)
$$
 рыухсЁ√ ╦ш Ёрэур 3 эрф хтъышфют√ь ъюы№Ўюь. ▌ЄюЄ ртЄюьюЁЇшчь  ты хЄё  рэрыюуюь ртЄюьюЁЇшчьр └эшър ётюсюфэ√ї рёёюЎшрЄштэ√ї рыухсЁ \cite{Cohn72}.

┬ ъюэЎх ётюхщ ЁрсюЄ√ \cite{Shevelin2} ╠.└. ╪хтхышэ юЄьхЄшы, ўЄю шч чрьхўрэшщ ЁхЎхэчхэЄр ёыхфєхЄ ёыхфє■∙шщ Ёхчєы№ЄрЄ: <<├Ёєяяр ртЄюьюЁЇшчьют ётюсюфэющ рыухсЁ√ ╦ш юЄ ЄЁхї яхЁхьхээ√ї  ты хЄё  ётюсюфэ√ь яЁюшчтхфхэшхь ё юс·хфшэхээющ яюфуЁєяяющ. ╥юўэ√щ тшф ёюьэюцшЄхыхщ ш юс·хфшэхээющ яюфуЁєяя√ ьюцэю шчтыхў№ шч т√ъырфюъ яєэъЄр 2.3>>. ═хфртэю ╨. ═рєЁ√чсрхт ш ╙. ╙ьшЁсрхт рэюэёшЁютрыш \cite{NU18}, ўЄю уЁєяяр ртЄюьюЁЇшчьют ётюсюфэющ рыухсЁх ╦ш $L_3(K)$ юЄ ЄЁхї яхЁхьхээ√ї эрф яЁюшчтюы№э√ь яюыхь  $K$ яЁхфёЄрты хЄё  т тшфх
$$
Aut(L_3(K)) \cong GL_3(K)*_{H(K)} B_3(K),
$$
уфх $GL_3(K)$ -- яюфуЁєяяр ышэхщэ√ї ртЄюьюЁЇшчьют ш $B_3(K)$ -- яюфуЁєяяр ртЄюьюЁЇшчьют тшфр
$$
\mu=(\alpha_1 x_1+h(x_2, x_3), \beta_2 x_2+ \beta_3 x_3, \gamma_2 x_2+\gamma_3 x_3), \;\;
$$
ш $H(K)=GL_3(K) \cap B_3(K)$.
 ┬ $\S$3 ь√ ЁрёяЁюёЄЁрэшыш ¤ЄюЄ Ёхчєы№ЄрЄ эр уЁєяяє Ёєўэ√ї ртЄюьюЁЇшчьют ётюсюфэ√ї рыухсЁ ╦ш (ш ётюсюфэ√ї рэЄшъюььєЄрЄштэ√ї рыухсЁ) юЄ ЄЁхї яхЁхьхээ√ї эрф яЁюшчтюы№эющ юсырёЄ№■ ЎхыюёЄэюёЄш $\Phi$.

┬ $\S$4 сыруюфрЁ  ¤Єюьє Ёхчєы№ЄрЄє эрь єфрыюё№ єёЄрэютшЄ№ ёюъЁрЄшьюёЄ№ Ёєўэ√ї ртЄюьюЁЇшчьют. ┬ $\S$5 яЁхфёЄртыхэю фюърчрЄхы№ёЄтю Єюую, ўЄю ртЄюьюЁЇшчь $\delta$ рыухсЁ√ $L_3(\Phi)$ эрф хтъышфют√ь ъюы№Ўюь $\Phi$  ты хЄё  фшъшь. ┬ $\S$2 яЁштхфхэ√ эхъюЄюЁ√х эхюсїюфшь√х юяЁхфхыхэш  ш ЇръЄ√.

\section{╬яЁхфхыхэш  ш яЁхфтрЁшЄхы№э√х ЇръЄ√}

\hspace*{\parindent}

╧єёЄ№ $\Phi$ - яЁюшчтюы№эр  юсырёЄ№ ЎхыюёЄэюёЄш. ╠эюцхёЄтю тёхї юсЁрЄшь√ї ¤ыхьхэЄют $\Phi$ юсючэрўшь ўхЁхч $\Phi^*$.
╬сючэрўшь ўхЁхч $L_{\Phi}\left\langle x_1, x_2, \ldots,x_n\right\rangle$ ётюсюфэє■ рыухсЁє ╦ш ё ьэюцхёЄтюь ётюсюфэ√ї яюЁюцфр■∙шї $X=\left\{x_1, x_2, \ldots, x_n\right\}$ эрф $\Phi$.
╧юыюцшь $x_1>x_2>\ldots>x_n$. ╬сючэрўшь ўхЁхч $X^*$ ьэюцхёЄтю тёхї эхрёёюЎшрЄштэ√ї ёыют т рыЇртшЄх $X$. {\em ─ышэющ
эхрёёюЎшрЄштэюую ёыютр} $v \in X^*$ эрч√трхЄё  ўшёыю тїюцфхэшщ т эхую ¤ыхьхэЄют ьэюцхёЄтр $X$, юэр юсючэрўрхЄё  ўхЁхч $d(v)$. ╩рцфюх эхрёёюЎшрЄштэюх ёыютю $v$ фышэ√ $\geq 2$ шьххЄ хфшэёЄтхээюх яЁхфёЄртыхэшх т тшфх $v=v_1v_2$, уфх $d(v_1), d(v_2)<d(v)$ \cite{ZSSSh}.

┬ ЁрсюЄх └.╚. ╪шЁ°ютр \cite{Shirshov53} фюърчрэю, ўЄю $L_{\Phi}\left\langle x_1, x_2, \ldots,x_n\right\rangle$  ты хЄё  ётюсюфэ√ь $\Phi$-ьюфєыхь эр ьэюцхёЄтх яЁртшы№э√ї ёыют. ▌ЄюЄ Ёхчєы№ЄрЄ яючтюы хЄ ЁрёёьрЄЁштрЄ№ $L_{\Phi}\left\langle x_1, x_2, \ldots,x_n\right\rangle$ ъръ $\Phi$-яюфьюфєы№ ётюсюфэющ рыухсЁ√ ╦ш $L_{K}\left\langle x_1, x_2, \ldots,x_n\right\rangle$ эрф яюыхь ўрёЄэ√ї $K=Q(\Phi)$ ъюы№Ўр $\Phi$.

┴єфхь уютюЁшЄ№, ўЄю эхрёёюЎшрЄштэюх ёыютю $v\in X^*$ шьххЄ Єшя $(m_1, m_2,\ldots, m_n)\in \mathbb{Z}^n_+$, уфх $\mathbb{Z}_+$ - ьэюцхёЄтю эхюЄЁшЎрЄхы№э√ї Ўхы√ї ўшёхы, хёыш эхрёёюЎшрЄштэюх ёыютю $v$ ёюфхЁцшЄ $x_i$ Ёютэю $m_i$ Ёрч. ╫шёыю $m_i$ сєфхь эрч√трЄ№ ёЄхяхэ№■ ьюэюьр $v$ яю $x_i$ ш сєфхь юсючэрўрЄ№ ъръ $\deg_{x_i} v$. ╧юыюцшь $\deg(v)=\deg_{x_1}(v)+\deg_{x_2}(v)+\cdots +\deg_{x_n}(v)$.

┴юыхх Єюую, фы  ы■сюую $w=(w_1, w_2, \ldots, w_n) \in \mathbb{Z}^n$ юяЁхфхышь $w$-ёЄхяхээє■ ЇєэъЎш■ $\deg_w $ яюырур 
$$\deg_w (x_1)= w_1, \deg_w(x_2)= w_2,\ldots ,\deg_w(x_n)=w_n.$$
 ┼ёыш $u=u_1 u_2$, уфх $d(u_1), d(u_2)< d(u)$, Єю яюыюцшь
$$
\deg_w(u)=\deg_w(u_1)+\deg_w(u_2).
$$

╥ръшь юсЁрчюь, ы■сюх $w \in \mathbb{Z}^n$ юяЁхфхы хЄ уЁрфєшЁютъє
\bes
L_{\Phi}\left\langle x_1, x_2, \ldots,x_n\right\rangle=\bigoplus_{k\in \mathbb{Z}} L_k,
\ees
 уфх $L_k$ - ышэхщэр  юсюыюўър ёыют ёЄхяхэш $k$ яю $\deg_w$. ╩рцф√щ  эхэєыхтющ ¤ыхьхэЄ $f\in L_{\Phi}\left\langle x_1, x_2, \ldots,x_n\right\rangle$ юфэючэрўэю яЁхфёЄрты хЄё  т тшфх
\bes
f=f_{k_1}+f_{k_2}+\ldots+f_{k_t}, \ \ k_1<k_2<\ldots<k_t, \ \  0\neq f_{k_j}\in L_{k_j}.
\ees
▌ыхьхэЄ $f_{k_t}$ эрч√трхЄё  {\em ёЄрЁ°хщ юфэюЁюфэющ ўрёЄ№■} ¤ыхьхэЄр $f$ яю юЄэю°хэш■ ъ ёЄхяхэш $\deg_w$.
┼ёыш $w=(1, 1, \ldots, 1) \in \mathbb{Z}^n$, Єю $\deg_w$ ёютярфрхЄ ё $\deg$. ╫хЁхч $\overline{f}$ сєфхь юсючэрўрЄ№ ёЄрЁ°є■ юфэюЁюфэє■ ўрёЄ№ $f$ яю юЄэю°хэш■ ъ ЇєэъЎшш ёЄхяхэш $\deg$.

─ы  ¤ыхьхэЄют $f_1, f_2, \ldots, f_r \in L_{\Phi}\left\langle x_1, x_2, \ldots,x_n\right\rangle$ ўхЁхч $\left\langle f_1, f_2, \ldots, f_r\right\rangle$ сєфхь юсючэрўрЄ№ $\Phi$-яюфрыухсЁє яюЁюцфхээє■ ¤Єшьш ¤ыхьхэЄрьш. ┼ёыш $\left\langle f_1, f_2, \ldots, f_r \right\rangle$ -- ётюсюфэр  рыухсЁр ╦ш юЄ яюЁюцфр■∙шї $f_1, f_2, \ldots, f_r$, Єю ьэюцхёЄтю ¤ыхьхэЄют $f_1, f_2, \ldots, f_r$ эрч√трхЄё  ╦ш-эхчртшёшь√ь (╦ш-ётюсюфэ√ь).

═ряюьэшь, ўЄю ьэюуююсЁрчшх рыухсЁ ╦ш эрф яюыхь $K$  ты хЄё  эшы№ёхэют√ь, Є.х., хёыш ¤ыхьхэЄ√ $f_1, f_2, \ldots, f_r $ рыухсЁ√ $ L_{K}\left\langle x_1, x_2, \ldots,x_n\right\rangle$ -- ╦ш-чртшёшь√, Єю эрщфхЄё  $f_i$ Єръющ, ўЄю $\overline{f_i}\in \left\langle \overline{f_1}, \ldots, \overline{f_{i-1}}, \overline{f_{i+1}}, \ldots, \overline{f_r} \right\rangle$.

─ц. ╦хтшэ \cite{Lewin} яюърчры, ўЄю яЁюшчтюы№эюх юфэюЁюфэюх ьэюуююсЁрчшх рыухсЁ  ты хЄё  °ЁрщхЁют√ь Єюуфр ш Єюы№ъю Єюуфр, ъюуфр юэю  ты хЄё  эшы№ёхэют√ь.

╤ыхфє■∙шщ ЇръЄ  ты хЄё  юс∙хшчтхёЄэ√ь \cite{Shirshov53}.

\begin{corollary}\label{c1}
▌ыхьхэЄ√ $f_1, f_2$ рыухсЁ√ $ L_{\Phi}\left\langle x_1, x_2, \ldots,x_n\right\rangle$ -- ╦ш-чртшёшь√ Єюуфр ш Єюы№ъю Єюуфр, ъюуфр $f_1, f_2$ ышэхщэю чртшёшь√ эрф $\Phi$.
\end{corollary}
\Proof
╧єёЄ№ $f_1, f_2$ -- эхэєыхт√х ╦ш-чртшёшь√х ¤ыхьхэЄ√ рыухсЁ√ $L_{\Phi}\left\langle x_1, x_2, \ldots,x_n\right\rangle$. ╥юуфр юэш юёЄр■Єё  ╦ш-чртшёшь√ьш т рыухсЁх
╦ш $L_{K}\left\langle x_1, x_2, \ldots,x_n\right\rangle$ эрф яюыхь  ўрёЄэ√ї $K=Q(\Phi)$ ъюы№Ўр $\Phi$.
┬ ёшыє эшы№ёхэютюёЄш ьэюуююсЁрчш  рыухсЁ ╦ш эрф яюыхь, шьххь $\overline{f_1} \in \left\langle \overline{f_2}\right\rangle_{K}$ шыш $\overline{f_2} \in \left\langle \overline{f_1}\right\rangle_{K}$, уфх $\left\langle \overline{f}\right\rangle_{K}$ юсючэрўрхЄ $K$-яюфрыухсЁє яюЁюцфхээє■ ¤ыхьхэЄюь $f$.

╥ръ ъръ $\left\langle f\right\rangle = K f$ фы  ы■сюую $f$ т рыухсЁх ╦ш, Єю, юЄё■фр яюыєўрхь $\overline{f_1}= \alpha \overline{f_2}$ шыш $\overline{f_2}= \alpha \overline{f_1}$, уфх $0\neq \alpha \in K$. ─юяєёЄшь $\overline{f_1}= \alpha \overline{f_2}$.
╨рёёьюЄЁшь ¤ыхьхэЄ√ $f_1- \alpha f_2, f_2$.  ╥юуфр
$$
\deg(\overline{f_1- \alpha f_2}) < \deg(\overline{f_2}).
$$
╥ръ ъръ ¤ыхьхэЄ√ $f_1- \alpha f_2, f_2$ ёэютр  ты ■Єё  ╦ш-чртшёшь√ьш, Єю, ъръ ш т√°х, шьххь
$
\overline{f_2} \in \left\langle \overline{f_1- \alpha f_2}\right\rangle =K (\overline{f_1- \alpha f_2})
$
шыш
$
\overline{f_1- \alpha f_2} \in \left\langle \overline{f_2}\right\rangle =K \overline{f_2}
$.
▌Єю тючьюцэю Єюы№ъю яЁш $f_1- \alpha f_2=0$.$\Box$

╙ЄтхЁцфхэшх ёыхфє■∙хщ ыхьь√ Єръцх їюЁю°ю шчтхёЄэю \cite{Shirshov84}.
\begin{lemma} \label{l1}
╧єёЄ№ $f_1, f_2, \ldots,f_r \in L_{\Phi}\left\langle x_1, x_2, \ldots,x_n\right\rangle$ ш $H \in \left\langle f_1, f_2, \ldots, f_r\right\rangle$. ╥юуфр хёыш $\overline{f_1}, \overline{f_2}, \ldots, \overline{f_r}$ -- ╦ш-эхчртшёшь√, Єю $\overline{H} \in \left\langle \overline{f_1}, \overline{f_2}, \ldots, \overline{f_r} \right\rangle$.
\end{lemma}
\Proof
╧юыюцшь $\deg (f_i)=n_i$, уфх $1\leq i\leq r$. ╧єёЄ№ $0\neq h=h(z_1, z_2, \ldots, z_r) \in L_{\Phi}\left\langle z_1, z_2, \ldots, z_r\right\rangle$ Єръюх, ўЄю $H=h(f_1, f_2, \ldots, f_r)$.
╧юыюцшь $w=(n_1, n_2, \ldots, n_r)\in \mathbb{Z}^r$ ш ЁрёёьюЄЁшь т рыухсЁх $L_{\Phi}\left\langle z_1, z_2, \ldots, z_r\right\rangle$ ЇєэъЎш■ ёЄхяхэш $\deg_w$.
╥юуфр $h= h' + \tilde{h}$, уфх $\tilde{h}$ - ёЄрЁ°р  юфэюЁюфэр  ўрёЄ№ $h$ юЄэюёшЄхы№эю $\deg_w$ ш $\deg_w (h') < \deg_w (\tilde{h})$. ╧єёЄ№ $\deg_w (h)=k$.
╟рьхЄшь, ўЄю $f_i=f'_i + \overline{f_i}$ фы  тёхї $i$, уфх $\overline{f_i}$ - ёЄрЁ°р  юфэюЁюфэр  ўрёЄ№ $f_i$ яю $\deg$. ╥юуфр
\begin{multline*}
H=h(f_1, f_2, \ldots, f_r)= h'(f_1, f_2, \ldots, f_r)+\tilde{h}(f_1, f_2, \ldots, f_r)=\\
=h'(f_1, f_2, \ldots, f_r)+\tilde{h}(f'_1+\overline{f_1}, f'_2+\overline{f_2}, \ldots, f'_r+\overline{f_r})=g'+\tilde{h}(\overline{f_1}, \overline{f_2}, \ldots, \overline{f_r}),
\end{multline*}
уфх $\deg (g') < k$. ╥ръ ъръ $\overline{f_1}, \overline{f_2}, \ldots, \overline{f_r}$ -- ╦ш-эхчртшёшь√х, Єю $\tilde{h}(\overline{f_1}, \overline{f_2}, \ldots, \overline{f_r})$ эх Ёртхэ эєы■ ш шьххЄ ёЄхяхэ№ $k$ т ёшыє т√сюЁр $w$.
╤ыхфютрЄхы№эю, $\overline{H}=\tilde{h}(\overline{f_1}, \overline{f_2}, \ldots, \overline{f_r}) \in \left\langle \overline{f_1}, \overline{f_2}, \ldots, \overline{f_r} \right\rangle$.
$\Box$

\begin{corollary} \label{c2}
╧єёЄ№ $f_1, f_2$ -- ╦ш-эхчртшёшь√х ¤ыхьхэЄ√ рыухсЁ√ $L(\Phi)=L_{\Phi}\left\langle x_1, x_2,\ldots,x_n\right\rangle$. ╥юуфр фы  ы■сюую  $0 \neq g(x_1, x_2) \in [L(\Phi), L(\Phi)]$ шьххЄ ьхёЄю эхЁртхэёЄтю $$\deg (g(f_{1}, f_{2})) > \max \{ \deg (f_{1}), \deg (f_{2})\}.$$
\end{corollary}
\Proof ╟рьхэ   $\Phi$ эр яюых ўрёЄэ√ї $K=Q(\Phi)$, ьюцэю ёўшЄрЄ№, ўЄю $\Phi$  ты хЄё  яюыхь.  ┼ёыш $\overline{f_1}, \overline{f_2}$ -- ышэхщэю эхчртшёшь√, Єю юэш  ты ■Єё  Єръцх ╦ш-эхчртшёшь√ьш яю ёыхфёЄтш■ \ref{c1}. ╧ю ыхььх \ref{l1} шьххь
$\overline{g(f_{1}, f_{2})} \in \left\langle \overline{f_1},\overline{f_2}\right\rangle$. ┴юыхх Єюую, ърцф√щ ьюэюь $g(x_1, x_2)$ ёюфхЁцшЄ $x_1$ ш $x_2$, Єръ ъръ $0 \neq g(x_1, x_2) \in [L(\Phi), L(\Phi)]$. ╧ютЄюЁ   Ёрёёєцфхэш  яЁштхфхээ√х т фюърчрЄхы№ёЄтх ыхьь√ \ref{l1}, яюыєўрхь
$$\deg (g(f_{1}, f_{2})) \geq \deg(f_{1})+\deg(f_{2}).$$

╧єёЄ№ $\overline{f_1} = \beta \overline{f_2}$, $0\neq \beta \in \Phi$. ╧юыюцшь $f_1'=f_1 - \beta f_2$. ╚ьххь $\left\langle f_1, f_2 \right\rangle = \left\langle f_1', f_2\right\rangle$ ш, юўхтшфэю, эрщфхЄё  $h \in [L(\Phi), L(\Phi)]$ Єръющ, ўЄю $g(f_1, f_2) = h(f_1', f_2)$. ╥ръ ъръ $\deg f_1' <\deg f_2$, Єю ¤ыхьхэЄ√ $\overline{f_1'}$ ш $\overline{f_2}$ -- ышэхщэю эхчртшёшь√. ╧Ёшьхэ   Ёрёёєцфхэш  яЁхф√фє∙хую рсчрЎр, шьххь
$$\deg (g(f_{1}, f_{2})) \geq \deg(f_{1}')+\deg(f_{2}).$$
╙ЄтхЁцфхэшх ёыхфёЄтш  т√яюыэ ■Єё  т юсюшї ёыєўр ї.
 $\Box$

\section{╧ЁхфёЄртыхэшх уЁєяя√ ртЄюьюЁЇшчьют рыухсЁ√ $L_{\Phi}\left\langle x_1, x_2, x_3\right\rangle$}

\hspace*{\parindent}

╧єёЄ№ $\mathrm{Aut}(L(\Phi))$ --- уЁєяяр тёхї ртЄюьюЁЇшчьют ётюсюфэющ рыухсЁ√ ╦ш $L_{\Phi}\left\langle x_1, x_2, \ldots,x_n\right\rangle$ эрф яЁюшчтюы№эющ юсырёЄ№■ ЎхыюёЄэюёЄш $\Phi$.
╫хЁхч ${\phi=(f_1,f_2, \ldots, f_n)}$ юсючэрўшь ртЄюьюЁЇшчь рыухсЁ√ $L_{\Phi}\left\langle x_1, x_2, \ldots,x_n\right\rangle$ Єръющ, ўЄю $\phi(x_i)=f_i$, уфх $1\leq i\leq n$.

└тЄюьюЁЇшчь $\delta$ рыухсЁ√ $L_{\Phi}\left\langle x_1, x_2, \ldots,x_n\right\rangle$ эрч√трхЄё  \textit{¤ыхьхэЄрЁэ√ь}, хёыш эрщфхЄё  $i$ Єръюх, ўЄю $\delta(x_i)=\alpha x_i+ f$ ш $\delta(x_j ) = x_j$ фы  тёхї $j\neq i$, уфх
 $\alpha\in \Phi^*, f \in L_{\Phi}\left\langle x_1, x_2, \ldots, x_{i-1}, x_{i+1},\ldots, x_n \right\rangle$. ╧юфуЁєяяр $T(L(\Phi))$ уЁєяя√ $\mathrm{Aut}(L(\Phi))$ яюЁюцфхээр  тёхьш ¤ыхьхэЄрЁэ√ьш ртЄюьюЁЇшчьрьш эрч√трхЄё  \textit{яюфуЁєяяющ Ёєўэ√ї ртЄюьюЁЇшчьют}.
└тЄюьюЁЇшчь $\varphi\in \mathrm{Aut}(L(\Phi))$ эрч√трхЄё  {\em Ёєўэ√ь}, хёыш $\varphi\in T(L(\Phi))$, шэрўх $\varphi $ эрч√трхЄё  {\em фшъшь.}

─рыхх тё■фє сєфхь ЁрёёьрЄЁштрЄ№ ЄЁхїяюЁюцфхээє■ ётюсюфэє■ рыухсЁє ╦ш $L_3(\Phi)=L_{\Phi}\left\langle x_1, x_2, x_3\right\rangle$.
╧єёЄ№ $GL_3(\Phi)$ -- яюфуЁєяяр ышэхщэ√ї ртЄюьюЁЇшчьют уЁєяя√ $\mathrm{Aut}(L_3(\Phi))$.
╫хЁхч $B_3(\Phi)$ юсючэрўшь яюфуЁєяяє уЁєяя√ ртЄюьюЁЇшчьют рыухсЁ√ $L_3(\Phi)$  тшфр
$$
\mu=(\alpha_1 x_1+h(x_2, x_3), \beta_2 x_2+ \beta_3 x_3, \gamma_2 x_2+\gamma_3 x_3),
$$
уфх $\alpha_1, \beta_2, \beta_3, \gamma_2, \gamma_3 \in \Phi, \; h \in L_{\Phi}\left\langle x_2, x_3\right\rangle$.
╬ўхтшфэю, $\mu$  ты хЄё  ртЄюьюЁЇшчьюь Єюуфр ш Єюы№ъю Єюуфр, ъюуфр $\alpha_1 \in \Phi^*$ ш
$$\begin{vmatrix}
\beta_2 & \beta_3\\
\gamma_2 & \gamma_3\\
\end{vmatrix} \in \Phi^*.$$

\begin{lemma} \label{l2}
CшёЄхьр ¤ыхьхэЄют
$$
B_0=\left\{\tau=(x_1+h(x_2,x_3), x_2, x_3)| h(x_2, x_3)\in \left[L_3(\Phi), L_3(\Phi) \right]\right\}
$$
  ты хЄё  ёшёЄхьющ яЁхфёЄртшЄхыхщ ыхт√ї ёьхцэ√ї ъырёёют уЁєяя√ $B_3(\Phi)$ яю яюфуЁєяях $H(\Phi)={GL_3(\Phi) \cap B_3(\Phi)}$.
\end{lemma}

\Proof ╦■сющ ¤ыхьхэЄ $\mu \in B_3(\Phi)$ шьххЄ тшф
\[
\mu = (\alpha_1 x_1 + \alpha_2 x_2 + \alpha_3 x_3 + g(x_2, x_3), \beta_2 x_2+ \beta_3 x_3, \gamma_2 x_2+\gamma_3 x_3),
\]
 уфх $\alpha_1 \in \Phi^*$ ш $g(x_2, x_3) \in \left[L_3(\Phi), L_3(\Phi) \right]$. ╥юуфр $\mu = \tau \circ \gamma$, уфх
\[
\tau=(x_1 + \alpha_1^{-1} g(x_2, x_3), x_2, x_3), \  \gamma=(\alpha_1 x_1 + \alpha_2 x_2 + \alpha_3 x_3, \beta_2 x_2+ \beta_3 x_3, \gamma_2 x_2+\gamma_3 x_3).
\]
╬ўхтшфэю, $\tau \in B_0, \gamma \in H(\Phi)$.

╧єёЄ№ ЄхяхЁ№ $\tau, \tau' \in B_0$ ш $\tau \neq \tau'$. ┼ёыш  $\tau=(x_1+h(x_2,x_3), x_2, x_3)$ ш $\tau'=(x_1+h'(x_2,x_3), x_2, x_3)$, Єю
\begin{eqnarray*}
\tau^{-1} \circ \tau'  = (x_1-h(x_2, x_3), x_2, x_3) \circ (x_1+h'(x_2, x_3), x_2, x_3) = \\
=(x_1-h(x_2, x_3)+h'(x_2, x_3), x_2, x_3) \notin H(\Phi),
\end{eqnarray*}
 яюёъюы№ъє $h(x_2, x_3) \neq h'(x_2, x_3)$. $\Box$

╧єёЄ№ $A_0$ --- яЁюшчтюы№эр  ЇшъёшЁютрээр  ёшёЄхьр яЁхфёЄртшЄхыхщ (ёюфхЁцр∙р  хфшэшЎє) ыхт√ї ёьхцэ√ї ъырёёют уЁєяя√ $GL_3(\Phi)$ яю яюфуЁєяях ${H(\Phi)=GL_3(\Phi) \cap B_3(\Phi)}$.

\begin{lemma} \label{l3}
╦■сющ Ёєўэющ ртЄюьюЁЇшчь $\phi \in T(L_3(\Phi))$ яЁхфёЄртшь
т тшфх
\begin{equation}\label{f1}
\phi= \sigma_1 \circ \tau_1 \circ \sigma_2 \circ \tau_2 \circ \cdots \circ \sigma_k \circ \tau_k \circ \lambda,
\end{equation}
уфх $\sigma_i \in A_0, \sigma_2, \ldots, \sigma_k \neq id, \tau_i \in B_0, \tau_1, \ldots, \tau_{k}\neq id, \lambda \in GL_3(\Phi)$.
\end{lemma}
\Proof
╦■сющ Ёєўэющ ртЄюьюЁЇшчь $\phi$ яЁхфёЄрты хЄё  т тшфх
$$
\phi= \delta_1 \circ \delta_2 \circ \ldots \circ \delta_n,
$$
уфх $\delta_1, \delta_2, \ldots, \delta_n$ - ¤ыхьхэЄрЁэ√х ртЄюьюЁЇшчь√.
╦■сющ ¤ыхьхэЄрЁэ√щ ртЄюьюЁЇшчь шьххЄ тшф
$$
\lambda_1 \circ \tau \circ \lambda_2,
$$
уфх $\tau \in B_0, \lambda_1, \lambda_2 \in GL_3(\Phi)$.

╤ыхфютрЄхы№эю, шьххь
$$
\phi= \lambda_1 \circ \tau_1 \circ \lambda_2 \circ \tau_2 \circ \ldots \circ \lambda_n \circ \tau_n \circ \lambda_{n+1}
$$
уфх $\lambda_i \in GL_3(\Phi), \tau_i \in B_0$.
╚эфєъЎшхщ яю $n$ фюърцхь, ўЄю $\phi$ яЁхфёЄрты хЄё  т тшфх \eqref{f1} ё $k\leq n$. ╚ьххь $\lambda_1=\sigma_1 \circ \gamma_1$, уфх $\sigma_1 \in A_0, \gamma_1 \in H(\Phi)$. ╥юуфр
$$
\lambda_1 \circ \tau_1=\sigma_1 \circ \gamma_1\circ \tau_1=\sigma_1 \circ \tau'_1\circ \gamma_1,
$$
уфх $\tau'_1=\gamma_1 \circ \tau_1\circ \gamma^{-1}_1 \in B_0$.
╤ыхфютрЄхы№эю,
$$
\phi= \sigma_1 \circ \tau'_1\circ (\gamma_1 \circ \lambda_2) \circ \tau_2 \circ \ldots \circ \lambda_n \circ \tau_n \circ \lambda_{n+1}.
$$
╧ю шэфєъЄштэюьє яЁхфяюыюцхэш■ яЁюшчтхфхэшх
$$
(\gamma_1 \circ \lambda_2) \circ \tau_2 \circ \ldots \circ \lambda_n \circ \tau_n \circ \lambda_{n+1}
$$
чряшё√трхЄё  т тшфх
$$
\sigma_2 \circ \tau'_2 \circ \sigma_3 \cdots \circ \sigma_k \circ \tau'_k \circ \sigma_{k+1} \circ \gamma, k\leq n.
$$
╤ыхфютрЄхы№эю,
$$
\phi=\sigma_1 \circ \tau'_1\circ\sigma_2 \circ \tau'_2 \circ \sigma_3 \cdots \circ \sigma_k \circ \tau'_k \circ \sigma_{k+1} \circ \gamma.
$$

┼ёыш $\sigma_2 \neq id$, Єю яюыєўхээюх яЁхфёЄртыхэшх шьххЄ тшф \eqref{f1}. ╧єёЄ№ $\sigma_2 = id$. ╥ръ ъръ $\tau'_1 \circ \tau'_2=\tau''_2 \in B_0$, Єю
$$
\phi=\sigma_1 \circ \tau'_1 \circ \tau'_2 \circ \sigma_3 \cdots \circ \sigma_k \circ \tau'_k \circ \sigma_{k+1} \circ \lambda=\sigma_1 \circ \tau''_2 \circ \sigma_3 \cdots \circ \sigma_k \circ \tau'_k \circ \sigma_{k+1} \circ \lambda.
$$
╧юёъюы№ъє $k-1<n$, Єю яю шэфєъЄштэюьє яЁхфяюыюцхэш■ $\phi$ чряшё√трхЄё  т тшфх \eqref{f1}.
$\Box$

┼ёыш $\phi=(f_1, f_2, f_3) \in T(\Phi)$,
Єю яюыюцшь
$$\deg (\phi) = \max \{\deg(f_1), \deg(f_2), deg(f_3)\}.$$

╨рёёьюЄЁшь яЁхфёЄртыхэшх  $\phi$ т тшфх \eqref{f1}. ─ы  ърцфюую $1 \leq i \leq k$ яюыюцшь
\begin{equation} \label{f2}
\phi_i = \sigma_1 \circ \tau_1 \circ \sigma_2 \circ \tau_2 \circ \ldots \circ \sigma_i \circ \tau_i = (\phi_i(x_1), \phi_i(x_2), \phi_i(x_3)).
\end{equation}

\begin{lemma} \label{l4}
─ы  ърцфюую ртЄюьюЁЇшчьр $\phi_i$ шч (\ref{f2}) шьххЄ ьхёЄю эхЁртхэёЄтю
\[
1<\deg (\phi_1) < \deg (\phi_2) < \ldots < \deg (\phi_k),
\]
\[
\deg (\phi_i(x_2)), \deg (\phi_i(x_3)) < \deg (\phi_i(x_1))=\deg (\phi_i).
\]
\end{lemma}

\Proof ─юърчрЄхы№ёЄтю фрээющ ыхьь√ яЁютхфхь шэфєъЎшхщ яю $k$. ╧єёЄ№ $k=1$ ш
$$
\tau_1 = (x_1+h_1(x_1, x_2), x_2, x_3),\;\;\; 0 \neq h_k(x_2, x_3) \in \left[L(\Phi),L(\Phi) \right].
$$
╥юуфр
$$
\phi_1 = \sigma_1 \circ \tau_1 = (\sigma_1(x_1)+h_1(\sigma_1(x_2), \sigma_1(x_3)), \sigma_1(x_2), \sigma_1(x_3)) = (\phi_1(x_1), \phi_1(x_2), \phi_1(x_3)).
$$
╥ръ ъръ $\sigma_1(x_1), \sigma_1(x_2), \sigma_1(x_3)$ -- ышэхщэ√х ш ышэхщэю эхчртшёшь√х, Єю
$$
\deg(h_1(\sigma_1(x_2), \sigma_1(x_3)))\geq 2.
$$
╤ыхфютрЄхы№эю, шьххь
$$
1 = \deg (\phi_1(x_2)) = \deg (\phi_1(x_3)) < \deg (\phi_1(x_1)) \ \text{ш} \ \deg (\phi_1) = \deg (\phi_1(x_1)).
$$
─юяєёЄшь, ўЄю єЄтхЁцфхэшх ыхьь√ тхЁэю фы  тёхї ртЄюьюЁЇшчьют $\phi_i$, уфх $\ 1 \leq i \leq k-1$. ╧юыюцшь
$$
\phi_{k-1}=(\phi_{k-1}(x_1), \phi_{k-1}(x_2), \phi_{k-1}(x_3))=(g,s,t).
$$
╧ю яЁхфяюыюцхэш■ шэфєъЎшш шьххь
$$
\deg (s), \deg (t) < \deg (g) \;\text{ш} \
\deg (\phi_{k-1}) =\deg (\phi_{k-1}(x_1))=\deg (g).
$$

╧єёЄ№ $$\sigma_{k}=(a_1x_1+b_1x_2+c_1x_3, a_2x_1+b_2x_2+c_2x_3, a_3x_1+b_3x_2+c_3x_3) \neq id.$$ ╥юуфр їюЄ  с√ юфшэ шч ¤ыхьхэЄют $a_2$ ш $a_3$ эх Ёртэ√ эєы■, Єръ ъръ $\sigma_{k}\notin H(\Phi)$. ╚ьххь
$$
\phi_{k-1}\circ \sigma_{k}=(a_1g+b_1s+c_1t, a_2g+b_2s+c_2t, a_3g+b_3s+c_3t).
$$
╤ыхфютрЄхы№эю, $\max\{\deg(a_2g+b_2s+c_2t), \deg(a_3g+b_3s+c_3t)\}=\deg (g)$.

╧єёЄ№
$$
\tau_k=(x_1+h_k(x_2, x_3), x_2, x_3), 0 \neq h_k(x_2, x_3) \in \left[L(\Phi),L(\Phi) \right].
$$
╥юуфр
\begin{multline*}
\phi_{k-1} \circ \sigma_{k} \circ \tau_{k}=
(a_1g+b_1s+c_1t+h_k(a_2g+b_2s+c_2t, a_3g+b_3s+c_3t),\\
 a_2g+b_2s+c_2t, a_3g+b_3s+c_3t)=(\phi_k(x_1), \phi_k(x_2), \phi_k(x_3))
\end{multline*}
╥ръ ъръ ъюьяюэхэЄ√
$$
\phi_k(x_2)=a_2g+b_2s+c_2t, \phi_k(x_3)=a_3g+b_3s+c_3t
$$
 ты ■Єё  ╦ш-эхчртшёшь√ьш, Єю яю ёыхфёЄтш■ 2 яюыєўрхь ёыхфє■∙хх эхЁртхэёЄтю
$$
\deg(h_k(a_2g+b_2s+c_2t, a_3g+b_3s+c_3t))>\max\{\deg(\phi_k(x_2)), \deg(\phi_k(x_3))\} =\deg(g).
$$
╤ыхфютрЄхы№эю, $\deg \phi_k(x_1) > \deg g$ ш $\deg \phi_{k} = \deg \phi_k(x_1) > \deg \phi_{k-1}$.
$\Box$

\begin{lemma} \label{l5}
╧ЁхфёЄртыхэшх \eqref{f1} ртЄюьюЁЇшчьр $\phi$  шч ыхьь√ \ref{l3}  ты хЄё  юфэючэрўэ√ь.
\end{lemma}
\Proof ─юёЄрЄюўэю яюърчрЄ№, ўЄю
$$\sigma_1\circ\tau_1\circ\sigma_2\circ\tau_2\circ\ldots\circ\sigma_k\circ \tau_k\circ\sigma_{k+1}\circ\lambda\neq id,$$
яЁш $k\geq 1$, $\sigma_i\in A_0,\:\sigma_2,\ldots,\sigma_k\neq id$, $\tau_i\in B_0,\:\tau_1,\ldots,\tau_k\neq id$, $\lambda\in GL_3(\Phi)$.

─юърцхь юЄ яЁюЄштэюую. ─юяєёЄшь
$$\sigma_1\circ\tau_1\circ\sigma_2\circ\tau_2\circ\ldots\circ\sigma_k\circ\tau_k\circ\sigma_{k+1}\circ\lambda= id.$$
╥юуфр шьххь
\begin{gather} \label{f3}
\tau_1\circ\sigma_2\circ\tau_2\circ\ldots\circ\sigma_k\circ\tau_k=\sigma_1^{-1}\circ\lambda^{-1}\circ\sigma_{k+1}^{-1}.
\end{gather}
╤юуырёэю ыхььх \ref{l4} ртЄюьюЁЇшчь
$$\phi=(\phi_k(x_1), \phi_k(x_2), \phi_k(x_3))=\tau_1\circ\sigma_2\circ\tau_2\circ\ldots\circ\sigma_k\circ\tau_k$$
шьххЄ ёЄхяхэ№
$$\deg(\phi)=\deg \phi_k(x_1)>1,$$
Є.х. ртЄюьюЁЇшчь т ыхтющ ўрёЄш ЁртхэёЄтр \eqref{f3} эх  ты хЄё  ышэхщэ√ь. ▌Єю яЁюЄштюЁхўшЄ ЁртхэёЄтє \eqref{f3}, Єръ ъръ
яЁртр  ўрёЄ№ ЁртхэёЄтр \eqref{f3}  ты хЄё  ышэхщэющ.
$\Box$

\begin{theorem} \label{T1}
╧єёЄ№ $\Phi$ - яЁюшчтюы№эр  юсырёЄ№ ЎхыюёЄэюёЄш. ├Ёєяяр Ёєўэ√ї ртЄюьюЁЇшчьют $T(L_3(\Phi))$ рыухсЁ√ $L_3(\Phi)$
 ты хЄё  ётюсюфэ√ь яЁюшчтхфхэшхь яюфуЁєяя√ ышэхщэ√ї ртЄюьюЁЇшчьют $GL_3(\Phi)$ ш яюфуЁєяя√ $B_3(\Phi)$ ё юс·хфшэхээющ яюфуЁєяяющ $H(\Phi) = GL_3(\Phi)\cap B_3(\Phi)$, Є.х.
$$
T(L_3(\Phi)) \cong GL_3(\Phi)*_{H(\Phi)} B_3(\Phi).
$$
\end{theorem}

\Proof
╧юёъюы№ъє $A_0$ ш $B_0$ -- ёшёЄхь√ ыхт√ї ёьхцэ√ї ъырёёют $GL_3(\Phi)$ ш $B_3(\Phi)$ яю яюфуЁєяях $H$, Єю
яю ыхььх \ref{l3} ы■сющ Ёєўэющ ртЄюьюЁЇшчь рыухсЁ√ $L_3(\Phi)$ яЁхфёЄрты хЄё  т тшфх \eqref{f1}. ╧ю ыхььх \ref{l5} Єръюх яЁхфёЄртыхэшх  ты хЄё  юфэючэрўэ√ь. ▌Єю ючэрўрхЄ \cite{MKS}, ўЄю
$$T_3(\Phi) \cong GL_3(\Phi) \ast_{H(\Phi)} B_3(\Phi). \ \ \ \ \Box $$

\section{╬ ёюъЁрЄшьюёЄш Ёєўэ√ї ртЄюьюЁЇшчьют}
\hspace*{\parindent}
╓хыюёЄэюх ъюы№Ўю $\Phi$, эх  ты ■∙ххё  яюыхь, эрч√трхЄё  {\em хтъышфют√ь} \cite{Vinberg}, хёыш ёє∙хёЄтєхЄ ЇєэъЎш 
$ \left|\cdot\right|:\Phi\backslash \left\{0\right\} \longrightarrow \mathbb{Z}_{+}$
(эрч√трхьр  эюЁьющ), єфютыхЄтюЁ ■∙р  ёыхфє■∙шь єёыютш ь:

E1) фы  ы■с√ї $a,b\in \Phi\backslash \left\{0\right\}$, $\left|ab\right|\geq \left|a\right|$, яЁшўхь ЁртхэёЄтю шьххЄ ьхёЄю Єюы№ъю Єюуфр, ъюуфр ¤ыхьхэЄ $b$ юсЁрЄшь;

E2) фы  ы■с√ї $a,b\in \Phi$, уфх $b \neq 0$, ёє∙хёЄтє■Є Єръшх $q, r\in \Phi$, ўЄю $a=bq + r$ ш ышсю $r=0$, ышсю $\left|r\right| < \left|b\right|$.

╧юыюцшь $e=|1| \in \mathbb{Z}_+$, уфх $1 \in \Phi$-хфшэшЎр ъюы№Ўр $\Phi$.
╚ьххь $|a|=e$ Єюуфр ш Єюы№ъю Єюуфр, ъюуфр $a\in \Phi^*$.

╬ёэютэ√ьш  яЁшьхЁрьш хтъышфют√ї ъюыхЎ  ты ■Єё  ъюы№Ўю $\mathbb{Z}$ Ўхы√ї ўшёхы ё рсёюы■Єэ√ь чэрўхэшхь Ўхы√ї ўшёхы ш ъюы№Ўю $K[x]$ ьэюуюўыхэют эрф яюыхь $K$ ёю ёЄхяхэ№■ ьэюуюўыхэют. ╤ыхфютрЄхы№эю, $e=1$ т ъюы№Ўх $\mathbb{Z}$ ш $e=0$ т ъюы№Ўх $K[x]$.

╧єёЄ№ $L_3(\Phi)=L_{\Phi}\left\langle x_1,x_2, x_3\right\rangle$ -- ётюсюфэр  рыухсЁр ╦ш юЄ ЄЁхї яхЁхьхээ√ї эрф хтъышфют√ь ъюы№Ўюь $\Phi$.  ┬тхфхь ышэхщэ√щ яюЁ фюъ $\geq$ эр ьэюцхёЄтх $X^*$ эхрёёюЎшрЄштэ√ї ёыют. ╧юыюцшь $x_1 > x_2>x_3$.
┼ёыш $u, v \in X^*$, Єю яюыюцшь  $u<v$, хёыш т√яюыэ хЄё  юфэю шч ёыхфє■∙шї єёыютшщ:

(i) $d(u)<d(v);$

(ii) ┼ёыш $d(u)=d(v)\geq2, u=u_1u_2$ ш $v=v_1v_2$, Єю $u_1<v_1$ шыш $u_1=v_1$ ш $u_2<v_2$.

╧ю ЄхюЁхьх └.╚. ╪шЁ°ютр ётюсюфэр  рыухсЁр ╦ш $L_{\Phi}\left\langle x_1, x_2,x_3\right\rangle$  ты хЄё  ётюсюфэ√ь $\Phi$-ьюфєыхь ё срчющ $V$, ёюёЄю ∙хщ шч яЁртшы№э√ї эхрёёюЎшрЄштэ√ї ёыют \cite{Shirshov53}.
╩рцф√щ эхэєыхтющ ¤ыхьхэЄ $f\in L_3(\Phi)$ чряшё√трхЄё  юфэючэрўэю т тшфх

$$
f=\alpha_1 v_1 +\alpha_2 v_2+\ldots+ \alpha_n v_n,\; 0\neq\alpha_i\in\Phi,  v_1>v_2>\ldots >v_n\in V
$$
╤ыютю $v_1$ эрч√трхЄё  ёЄрЁ°шь ёыютюь (ьюэюьюь) $f$, р $\alpha_1$ эрч√трхЄё  ёЄрЁ°шь ъю¤ЇЇшЎшхэЄюь $f$. ╬сючэрўшь шї ўхЁхч $lm(f)$ ш $lc(f)$, ёююЄтхЄёЄтхээю. ╫хЁхч $\widehat{f}=\alpha_1 v_1$ юсючэрўшь ёЄрЁ°шщ ўыхэ ¤ыхьхэЄр $f$.

╧єёЄ№ $f$ --- яЁюшчтюы№э√щ ¤ыхьхэЄ шч $L_3(\Phi)$. ╧юыюцшь
$$
D(f)=(lm(f),|lc(f)|)
$$
 ш эрчютхь $D(f)$ яюърчрЄхыхь ¤ыхьхэЄр $f$.

╧єёЄ№ $\varphi =(f_1, f_2, f_3)$ - ртЄюьюЁЇшчь рыухсЁ√ $L_3(\Phi)$.
╥юуфр яюърчрЄхыхь  ртЄюьюЁЇшчьр $\varphi$ эрчютхь $$D(\varphi)=(u, w, s, |lc(f_1)|+|lc(f_2)|+|lc(f_3)|)\in V^3\times Z_{+},$$
уфх
$\left\{u, w, s\right\}= \left\{lm(f_1),lm(f_2), lm(f_3)\right\}$  ш $u \geq w \geq s$.

╟рьхЄшь, ўЄю $D$ шэтрЁшрэЄэю юЄэюёшЄхы№эю яхЁхёЄрэютъш ъюьяюэхэЄ ртЄюьюЁЇшчьр.
╚ьххь
$$
D(id)=(x_1, x_2, x_3, 3e),
$$
уфх $id=(x_1, x_2, x_3)$ - ЄюцфхёЄтхээ√щ ртЄюьюЁЇшчь.
┴юыхх Єюую, $$D(f_1, f_2, f_3)=(x_1, x_2, x_3, 3e)$$ Єюуфр ш Єюы№ъю Єюуфр, ъюуфр ¤ыхьхэЄ√ $f_1, f_2, f_3$ ёютярфр■Є ё эхъюЄюЁющ яхЁхёЄрэютъющ ¤ыхьхэЄют
\begin{equation}\label{f4}
\alpha_{1} x_1 + \beta_{1} x_2 + \gamma_{1} x_3, \ \alpha_{2} x_2 + \beta_{2} x_3, \ \alpha_{3}x_3,
\end{equation}
уфх $\alpha_1, \alpha_2, \alpha_3 \in\Phi^*$ ш $ \beta_1, \beta_2, \gamma_1 \in\Phi$.

╬сючэрўшь ўхЁхч $\preceq$ ыхъёшъюуЁрЇшўхёъшщ яюЁ фюъ эр ьэюцхёЄтх $V^3\times Z_{+}$. ╟рьхЄшь, ўЄю ьэюцхёЄтю $V^3\times Z_{+}$ ышэхщэю єяюЁ фюўхэю юЄэюёшЄхы№эю $\preceq$.

┬тхфхь эхюсїюфшь√х яюэ Єш  ш юяЁхфхыхэш  ю ёюъЁрЄшьюёЄш ртЄюьюЁЇшчьют.

{\em ▌ыхьхэЄрЁэ√ь яЁхюсЁрчютрэшхь } ёшёЄхь√ ¤ыхьхэЄют $(f_1, f_2, f_3)$ эрч√трхЄё  чрьхэр юфэюую ¤ыхьхэЄр $f_i$ эр ¤ыхьхэЄ тшфр $\alpha f_i +g$, уфх $\alpha \in \Phi^*, g\in \left\langle \left\{f_j |j\ne i \right\} \right\rangle$.

╟ряшё№
$$
(f_1, f_2, f_3)\to (g_1, g_2, g_3)
$$
ючэрўрхЄ, ўЄю ёшёЄхьр ¤ыхьхэЄют $(g_1, g_2, g_3)$ яюыєўхэр шч ёшёЄхь√ ¤ыхьхэЄют $(f_1, f_2, f_3)$ юфэшь ¤ыхьхэЄрЁэ√ь яЁхюсЁрчютрэшхь.

┼ёыш $(f_1,f_2, f_3)$ - Ёєўэющ ртЄюьюЁЇшчь рыухсЁ√ $L_3(\Phi)$, Єю ёє∙хёЄтєхЄ яюёыхфютрЄхы№эюёЄ№ ¤ыхьхэЄрЁэ√ї яЁхюсЁрчютрэшщ тшфр

$$
(x_1, x_2, x_3)=(f^{(0)}_1,f^{(0)}_2, f^{(0)}_3)\rightarrow (f^{(1)}_1,f^{(1)}_2, f^{(1)}_2)\rightarrow \ldots \rightarrow (f^{(k)}_1, f^{(k)}_2, f^{(k)}_2)=(f_1, f_2, f_3).
$$

└тЄюьюЁЇшчь $\theta=(f_1, f_2, f_3)$ эрч√трхЄё  {\em ¤ыхьхэЄрЁэю $D$-ёюъЁрЄшь√ь} (шыш $\theta$ {\em фюяєёърхЄ ¤ыхьхэЄрЁэюх $D$-ёюъЁр∙хэшх}), хёыш ёє∙хёЄтєхЄ ртЄюьюЁЇшчь $\psi$ Єръющ, ўЄю ${\theta \to \psi}$ ш $D(\psi)\prec D(\theta)$.
┴єфхь уютюЁшЄ№, ўЄю {\em ртЄюьюЁЇшчь $\psi$  ты хЄё  $D$-¤ыхьхэЄрЁэ√ь ёюъЁр∙хэшхь ртЄюьюЁЇшчьр $\theta$.}

\begin{lemma}\label{l6}
╧єёЄ№ $\pi=(g_1, g_2, g_3)$ - ртЄюьюЁЇшчь рыухсЁ√ $L_3(\Phi)$. ┼ёыш ёЄрЁ°шх ўыхэ√ ¤ыхьхэЄют $g_i, g_j$ ышэхщэю чртшёшь√, уфх $1\leq i\neq j \leq 3$ , Єю ртЄюьюЁЇшчь $\pi$  ты хЄё  ¤ыхьхэЄрЁэю $D$-ёюъЁрЄшь√ь.
\end{lemma}
\Proof
┴хч юуЁрэшўхэш  юс∙эюёЄш яЁхфяюыюцшь, ўЄю $\widehat{g_1}$ ш $\widehat{g_2}$ ышэхщэю чртшёшь√, Єю хёЄ№ $lm(g_1) =lm(g_2)$, ш яєёЄ№ фы  юяЁхфхыхээюёЄш  $|lc(g_1)|\geq |lc(g_2)|$.
╧ю E2) ёє∙хёЄтє■Є  $q,r\in \Phi$ Єръшх, ўЄю ${lc(g_1)= lc(g_2) q+r}$ ш ышсю $r=0$, ышсю $|r|<|lc(g_2)|$.
╨рёёьюЄЁшь ¤ыхьхэЄрЁэюх яЁхюсЁрчютрэшх
$$\pi=(g_1, g_2, g_3)\to (g_1-q g_2, g_2, g_3)=\delta.$$
╚ьххь $D(g_1)\succ D(g_1-q g_2)$.
╤ыхфютрЄхы№эю, $D(\pi)\succ D(\delta)$ ш ртЄюьюЁЇшчь $\pi$  ты хЄё  ¤ыхьхэЄрЁэю $D$-ёюъЁрЄшь√ь.
$\Box$

╒рЁръЄхЁшчрЎш■ Ёєўэ√ї ртЄюьюЁЇшчьют рыухсЁ√ $L_3(\Phi)$ фрхЄ ёыхфє■∙хх яЁхфыюцхэшх.

\begin{proposition}\label{p1}
╧єёЄ№ $\phi=(f_1, f_2, f_3)$ - Ёєўэющ ртЄюьюЁЇшчь рыухсЁ√ $L_3(\Phi)$. ┼ёыш
$$D(\phi)\succ (x_1, x_2, x_3, 3e),$$
Єю ртЄюьюЁЇшчь $\phi$  ты хЄё  ¤ыхьхэЄрЁэю $D$-ёюъЁрЄшь√ь.
\end{proposition}

\Proof
╧ю ыхььх \ref{l5}
$\phi$ юфэючэрўэю чряшё√трхЄё  т тшфх \eqref{f1}.
╨рёёьюЄЁшь ёыєўрщ ъюуфр
$\lambda=id$.
╥юуфр
$$
\phi=\sigma_1 \circ \tau_1 \circ \sigma_2 \circ \tau_2 \circ \cdots \circ \sigma_k \circ \tau_k=(f_1, f_2, f_3).
$$
╧юыюцшь
$$
\psi=\sigma_1 \circ \tau_1 \circ \sigma_2 \circ \tau_2 \circ \cdots \circ \sigma_k=(g_1,g_2,g_3).
$$
┼ёыш $\tau_k=(x_1 + h_k(x_2,x_3), x_2, x_3)$, Єю
$$\phi=\psi \circ \tau_k=(g_1+h_k(g_2,g_3), g_2, g_3)=(f_1, f_2, f_3).$$
╧ю ыхььх \ref{l4} шьххь
$\deg (\psi)=\deg (\phi_{k-1})< \deg (\phi_k) =\deg (\phi)$. \\
╤ыхфютрЄхы№эю, $${D(\psi)\prec D(\phi).}$$
╧юёъюы№ъє $
\phi \to \psi,
$
Єю ртЄюьюЁЇшчь $\phi$  ты хЄё  ¤ыхьхэЄрЁэю $D$-ёюъЁрЄшь√ь.

─юяєёЄшь, ЄхяхЁ№, ўЄю
$$
\lambda=(a_1x_1+b_1x_2+c_1 x_3, a_2x_1+b_2x_2+c_2 x_3, a_3 x_1+b_3 x_2+c_3 x_3)\neq id,
$$
уфх шьххь
$$\begin{vmatrix}
a_1 & b_1 & c_1\\
a_2 & b_2 & c_2\\
a_3 & b_3 & c_3\\
\end{vmatrix} \in \Phi^*.$$
╧юыюцшь
$$
\omega= \sigma_1 \circ \tau_1 \circ \sigma_2 \circ \tau_2 \circ \cdots \circ \sigma_k \circ \tau_k=(g_1+h_k(g_2,g_3), g_2, g_3)=(w_1, w_2, w_3)
$$
╥юуфр, яю ыхььх \ref{l4} $\deg (w_1)> \max \left\{\deg (w_2), \deg (w_3)\right\}$.
╤ыхфютрЄхы№эю,
$$
\widehat{w_1}=\widehat{h_k(w_2, w_3)}.
$$
╚ьххь
$$
\phi=\omega \circ \lambda=(a_1 w_1+b_1 w_2 +c_1 w_3, a_2 w_1+b_2 w_2+c_2 w_3, a_3 w_1+b_3 w_2+c_3 w_3)=(f_1, f_2, f_3).
$$

┼ёыш їюЄ  с√ фтр ъю¤ЇЇшЎшхэЄр шч $a_1, a_2, a_3$ эхэєыхт√х, Єю ёЄрЁ°шх ёыютр їюЄ  с√ фтєї ъюьяюэхэЄ ртЄюьюЁЇшчьр $\phi$ ёютярфр■Є.
╥юуфр  $\phi$ - ¤ыхьхэЄрЁэю $D$-ёюъЁрЄшь яю ыхььх \ref{l6}.

┴хч юуЁрэшўхэшх юс∙эюёЄш ьюцэю ёўшЄрЄ№, ўЄю $a_1 = a_2 = 0,  a_3 \neq 0$.
╥юуфр шч
$$\begin{vmatrix}
0 & b_1 & c_1\\
0 & b_2 & c_2\\
a_3 & b_3 & c_3\\
\end{vmatrix} \in \Phi^*,$$
ёыхфєхЄ, ўЄю $a_3(b_1c_2-b_2c_1)\in \Phi^*$.
╤ыхфютрЄхы№эю $a_3\in \Phi^*$ ш $b_1c_2-b_2c_1\in \Phi^*$.

╚ьххь
$$
\left\langle b_1 w_2 + c_1 w_3, b_2 w_2 + c_2 w_3 \right\rangle = \left\langle w_2, w_3 \right\rangle.
$$
╨рёёьюЄЁшь ртЄюьюЁЇшчь
$$
\psi=(b_1 w_2 + c_1 w_3, b_2 w_2 + c_2 w_3, a^{-1}_3f_3 - h_k(w_2, w_3)),
$$
уфх $a^{-1}_3f_3 - h_k(w_2, w_3)= g_1 + a^{-1}_3(b_3 w_2 + c_3 w_3)$ ш $\deg(w_1) > \deg(g_1), \deg (w_2), \deg (w_3)$.
╧юёъюы№ъє
$\phi \to \psi$,
Єю ртЄюьюЁЇшчь $\phi$  ты хЄё  ¤ыхьхэЄрЁэю $D$-ёюъЁрЄшь√ь. $\Box$

═хяюёЁхфёЄтхээ√щ рэрышч фюърчрЄхы№ёЄтр яЁхфыюцхэш  \ref{p1} фрхЄ
\begin{corollary}\label{c3}
╧єёЄ№ $\pi=(g_1, g_2, g_3)$ -- Ёєўэющ ртЄюьюЁЇшчь рыухсЁ√ $L_3(\Phi)$ ш $lm(g_1)<lm(g_2)<lm(g_3)$. ╥юуфр эрщфхЄё  $h \in L_2(\Phi)$ Єръюх, ўЄю
$$
\overline{g_3}=\overline{h(g_1, g_2)}.
$$
\end{corollary}

\begin{corollary}\label{c4}
╨єўэ√х ш фшъшх ртЄюьюЁЇшчь√ рыухсЁ√ $L_3(\Phi)$ эрф ъюэёЄЁєъЄштэ√ь хтъышфют√ь ъюы№Ўюь $\Phi$ рыуюЁшЄьшўхёъш Ёрёяючэртрхь√.
\end{corollary}

\Proof
╧Ёютхфхь шэфєъЎш■ яю яюърчрЄхы■ $D(\phi)$ ртЄюьюЁЇшчьр ${\phi \in Aut(L_3(\Phi))}$.
┼ёыш $D(\phi)=(x_1, x_2, x_3, 3e)$, Єю $\phi$ шьххЄ тшф \eqref{f4}. ╤ыхфютрЄхы№эю, ртЄюьюЁЇшчь $\phi$  ты хЄё  ышэхщэ√ь. ╥ръ ъръ ы■ср  юсЁрЄшьр  ьрЄЁшЎр эрф хтъышфют√ь ъюы№Ўюь  ты хЄё  яЁюшчтхфхэшхь ¤ыхьхэЄрЁэ√ї ш фшруюэры№э√ї ьрЄЁшЎ \cite{Vinberg}, Єю ышэхщэ√х ртЄюьюЁЇшчь√  ты ■Єё  Ёєўэ√ьш.

┼ёыш ртЄюьюЁЇшчь $\phi$ эх  ты хЄё  ¤ыхьхэЄрЁэю $D$-ёюъЁрЄшь√ь, Єю $\phi$  ты хЄё  фшъшь яю яЁхфыюцхэш■ \ref{p1}. ╧єёЄ№ $\phi=(f_1, f_2, f_3)$ -- Ёєўэющ ртЄюьюЁЇшчь рыухсЁ√ $L_3(\Phi)$. ┼ёыш $lm(f_i), lm(f_j)$ ёютярфр■Є фы  эхъюЄюЁ√ї $1\leq i \neq j \leq 3$, Єю яю ыхььх \ref{l6} ртЄюьюЁЇшчь $\phi$ фюяєёърхЄ ышэхщэюх $D$-ёюъЁр∙хэшх. ╧ю¤Єюьє, схч юуЁрэшўхэш  юс∙эюёЄш ьюцэю ёўшЄрЄ№, ўЄю $lm(f_1) < lm(f_2) < lm(f_3)$.
╬Єё■фр ёыхфєхЄ, ўЄю $\overline{f_1}, \overline{f_2}, \overline{f_3}$ -- ышэхщэю эхчртшёшь√. ╧ю ёыхфёЄтш■ \ref{c3} эрщфхЄё  $h \in L_2(\Phi)$ Єръюх, ўЄю $\overline{f_3} = \overline{h(f_1, f_2)}$. ╟рьхЄшь, ўЄю $\overline{f_1}, \overline{f_2}$ -- ╦ш-эхчртшёшь√ яю ёыхфёЄтш■ \ref{c1}. ╥юуфр яю ыхььх \ref{l1} шьххь
$$
\overline{f_3} = \overline{h(f_1, f_2)} \in \left\langle \overline{f_1}, \overline{f_2}\right\rangle.
$$
╙ёыютшх $\overline{f_3} \in \left\langle \overline{f_1}, \overline{f_2}\right\rangle$ ыхуъю яЁютхЁ хЄё , Єръ ъръ $ \overline{f_1}, \overline{f_2}$ -- ╦ш-эхчртшёшь√ ш юфэюЁюфэ√.
┼ёыш $\overline{f_3} \in \left\langle \overline{f_1}, \overline{f_2}\right\rangle$, Єю ыхуъю яюёЄЁюшЄ№ $H\in L_2(\Phi)$ Єръюх, ўЄю $\overline{f_3} = H(\overline{f_1}, \overline{f_2})$.
╥юуфр ртЄюьюЁЇшчь $\psi=(f_1, f_2, f_3 - H(f_1, f_2))$  ты хЄё  ¤ыхьхэЄрЁэ√ь $D$-ёюъЁр∙хэшхь ртЄюьюЁЇшчьр $\phi$. ╚ьххь $D(\psi)\prec D(\phi)$. ╬ўхтшфэю, $\phi$  ты хЄё  Ёєўэ√ь Єюуфр ш Єюы№ъю Єюуфр, ъюуфр $\psi$  ты хЄё  Ёєўэ√ь ртЄюьюЁЇшчьюь. ╚эфєъЎш  яю $D$ чртхЁ°рхЄ фюърчрЄхы№ёЄтю.
$\Box$
\section{└эрыюу ртЄюьюЁЇшчьр └эшър}
\hspace*{\parindent}

╧єёЄ№ $\Phi$ --- хтъышфютю ъюы№Ўю ш $0\neq z\in \Phi\setminus \Phi^*$. ╧юыюцшь, ўЄю ${K=Q(\Phi)}$ --- яюых ўрёЄэ√ї ъюы№Ўр $\Phi$.
┬ ЁрсюЄх \cite{AU17} с√ыю яЁштхфхэю яюёЄЁюхэшх ртЄюьюЁЇшчьр ═рурЄ√ рыухсЁ√ $K[x_1,x_2]$ эрф яюыхь ўрёЄэ√ї $K$. ╥ръшь цх юсЁрчюь с√ы яюёЄЁюхэ ртЄюьюЁЇшчь ётюсюфэющ эхрёёюЎшрЄштэющ рыухсЁ√ Ёрэур 2 эрф хтъышфют√ь ъюы№Ўюь.
╬фэръю ьхЄюф яюёЄЁюхэш  ¤Єюую ртЄюьюЁЇшчьр эх яЁюїюфшЄ фы  рэЄшъюььєЄрЄштэ√ї рыухсЁ, Єръ ъръ ътрфЁрЄ ы■сюую ¤ыхьхэЄр т ¤Єющ рыухсЁх Ёртхэ 0. ╨рёёьюЄЁшь ёыхфє■∙є■ яюёыхфютрЄхы№эюёЄ№ яЁхюсЁрчютрэшщ ЄЁюхъ ¤ыхьхэЄют рыухсЁ√ $L_3(Q(\Phi))$:
\begin{multline*}
(x_1, x_2, x_3) \to (zx_1, x_2, x_3) \to (zx_1- [x_2, x_3], x_2, x_3) \to \\
 \to (zx_1- [x_2, x_3], x_2 + z(zx_1- [x_2,x_3]), x_3)\to \\
\to (zx_1- [x_2, x_3] +[x_2 + z(zx_1- [x_2, x_3]),x_3], x_2 + z(zx_1- [x_2,x_3]), x_3)=\\
=(zx_1+ z[zx_1- [x_2, x_3],x_3], x_2 + z(zx_1- [x_2, x_3]), x_3)=\\
=(z(x_1+ [zx_1- [x_2 x_3],x_3]), x_2 + z(zx_1- [x_2, x_3]), x_3)\to \\
\to (x_1+ [zx_1- [x_2, x_3],x_3], x_2 + z(zx_1- [x_2, x_3]), x_3).
\end{multline*}
╥ръшь юсЁрчюь ё яюью∙№■ ¤Єшї яЁхюсЁрчютрэшш, ъюЄюЁ√х  ты ■Єё  ¤ыхьхэЄрЁэ√ьш эрф $K$, эю эх эрф $\Phi$, яюыєўхэ ¤эфюьюЁЇшчь
$$
\delta=(x_1+ [zx_1- [x_2, x_3],x_3], x_2 + z(zx_1- [x_2, x_3]), x_3)
$$
рыухсЁ√ $L_3(\Phi)$.

\begin{lemma}\label{l7}
▌эфюьюЁЇшчь $\delta$ рыухсЁ√ ${L_3(\Phi)=L_\Phi\left\langle x_1,x_2, x_3\right\rangle}$  ты хЄё  ртЄюьюЁЇшчьюь.
\end{lemma}
\Proof
╧єёЄ№
$
\delta=(b_1,b_2, b_3).
$
 ─юёЄрЄюўэю яюърчрЄ№, ўЄю ¤ыхьхэЄ√ $b_1, b_2, b_3$ яюЁюцфр■Є тё■ рыухсЁє $L_\Phi\left\langle x_1, x_2, x_3\right\rangle$.
═хяюёЁхфёЄтхээ√х т√ўшёыхэш  фр■Є
\[zb_1-[b_2,b_3]=zx_1-[x_2,x_3] \in \left\langle b_1, b_2, b_3\right\rangle,\]
\[b_2-z(zb_1-[b_2,b_3])=x_2 \in \left\langle b_1, b_2, b_3\right\rangle,\]
\[b_1-[zb_1-[b_2,b_3],b_3]=x_1 \in \left\langle b_1, b_2, b_3\right\rangle,\]
Є.х., $x_1, x_2, x_3 \in \left\langle b_1, b_2, b_3\right\rangle$.
$\Box$

\begin{theorem}\label{t2}
 └тЄюьюЁЇшчь $\delta$ рыухсЁ√ $L_3(\Phi)=L_\Phi \left\langle x_1,x_2, x_3\right\rangle $
 ты хЄё  фшъшь.
\end{theorem}

\Proof
─юяєёЄшь, ўЄю ртЄюьюЁЇшчь
$$
\delta=(x_1+[zx_1-[x_2,x_3],x_3], x_2+z(zx_1-[x_2, x_3]), x_3)=(f_1, f_2, f_3)
$$
  ты хЄё  Ёєўэ√ь. ╚ьххь
$$
\overline{f_1}=[[x_2,x_3],x_3], \overline{f_2}=z[x_2,x_3], \overline{f_3}=x_3.
$$
╧ю ёыхфёЄтш■ \ref{c3} эрщфхЄё 
$ 0 \neq h \in L_\Phi \left\langle f_2, f_3\right\rangle$
Єръюх, ўЄю
 $\overline{f_1}=\overline{h(f_2,f_3)}$.
╥юуфр яю ыхььх \ref{l1} яюыєўрхь, ўЄю
$\overline{f_1}\in \left\langle \overline{f_2}, \overline{f_3}\right\rangle$.
┬ яюфрыухсЁх $\left\langle \overline{f_2}, \overline{f_3}\right\rangle$ Єюы№ъю яЁюшчтхфхэшх $[\overline{f_2}, \overline{f_3}]=z [[x_2, x_3], x_3]$ шьххЄ ёЄхяхэ№ 3
 яю $x_1, x_2, x_3$.
╥ръ ъръ $\overline{f_1}=[[x_2,x_3],x_3]$ ышэхщэю эх т√ЁрцрхЄё  ўхЁхч $[\overline{f_2}, \overline{f_3}]$, Єю $\overline{f_1}$ эх яЁшэрфыхцшЄ яюфрыухсЁх $\left\langle \overline{f_2}, \overline{f_3}\right\rangle$.
╧юыєўхээюх яЁюЄштюЁхўшх яюърч√трхЄ, ўЄю ртЄюьюЁЇшчь $\delta$ эх  ты хЄё  ¤ыхьхэЄрЁэю $D$-ёюъЁрЄшь√ь. ╧ю яЁхфыюцхэш■ 1 юэ эх  ты хЄё  Ёєўэ√ь.
$\Box$

╟рьхЄшь, ўЄю фюърчрЄхы№ёЄтр ЄхюЁхь 1 ш 2 яюыэюёЄ№■ яЁюїюф Є ш фы  ётюсюфэ√ї рэЄшъюььєЄрЄштэ√ї рыухсЁ. ┬ ўрёЄэюёЄш, ртЄюьюЁЇшчь  $\delta$
ётюсюфэющ рэЄшъюььєЄрЄштэющ рыухсЁ√ $AK_{\Phi}\left\langle x_1, x_2, x_3 \right\rangle$ эрф хтъышфют√ь ъюы№Ўюь $\Phi$ Єръцх  ты хЄё  фшъшь.


\begin{thebibliography}{99}

\bibitem{Umi07} U.U. Umirbaev,
{\it The Anick automorphism of free associative algebras}, J. Reine Angew. Math., {\bf 605} (2007), 165--178.

\bibitem{Jung} H.W.E. Jung,
{\it \"Uber ganze birationale Transformationen der Ebene}, J. Reine Angew. Math., {\bf 184} (1942), 161--174.

\bibitem{Kulk} W. van der Kulk,
{\it On polynomial rings in two variables}, Nieuw Arch. Wiskunde, (3){\bf 1} (1953), 33--41.

\bibitem{Shafarevich} I.R. Shafarevich, {\it On some infinite-dimensional groups}, Rend. Mat. e Appl. (5), 25:1--2 (1966), 208--212.

\bibitem{Wright} D. Wright,
{\it The amalgamated free product structure of $GL_2(k[x_1 ,...,x_n])$}, J. Pure Appl. Algebra, {\bf 12} (1978),
235--251

\bibitem{Makar-Limanov} ╦. ╠рърЁ-╦шьрэют,
{\it └тЄюьюЁЇшчь√ ётюсюфэющ рыухсЁ√ юЄ фтєї яюЁюцфр■∙шї}, ╘єэъЎшюэ. рэрышч ш хую яЁшы., {\bf 4}:3 (1970), 107--108; English translation: in Functional Anal. Appl., {\bf 4} (1970), 262--263.

\bibitem{Czerniakiewicz} A.G. Czerniakiewicz,
{\it Automorphisms of a free associative algebra of rank 2}, I, II, Trans. Amer. Math. Soc. {\bf 160} (1971), 393--401; {\bf 171} (1972), 309--315.


\bibitem{Poisson} L. Makar-Limanov, U. Turusbekova, U. Umirbaev,
{\it Automorphisms and derivations of free Poisson algebras in two variables}, J. Algebra, {\bf 322} (2009), 3318--3330.

\bibitem{right-symmetric} D. Kozybaev, L. Makar-Limanov, U. Umirbaev,
{\it The Freiheitssatz and the automorphisms of free right-symmetric algebras}, Asian-European J. Math., {\bf 1}:2 (2008), 243--254.

\bibitem{NK} ─.╒. ╩юч√срхт, └.╤. ═рєЁрчсхъютр,
{\it ╦шэхрЁшчрЎш  ртЄюьюЁЇшчьют ётюсюфэ√ї яЁртюёшььхЄЁшўэ√ї рыухсЁ
Ёрэур 2,} ╥хчшё√ фюъырфют ьхцфєэрЁюфэющ ъюэЇхЁхэЎшш ╠└╦▄╓┼┬╤╩╚┼ ╫╥┼═╚▀, (2018), 155--155.
http://www.math.nsc.ru/conference/malmeet/18/maltsev18.pdf



\bibitem{Umi04} I.P. Shestakov and U.U. Umirbaev,
{\it Tame and wild automorphisms of rings of polynomials in three variables}, J. Amer. Math. Soc. {\bf 17} (2004), 197--227.

\bibitem{Cohn} P.M. Cohn,
{\it Subalgebras of free associative algebras}, Proc. London Math. Soc., {\bf 56} (1964), 618--632.

\bibitem{Lewin} J. Lewin,
{\it On Schreier varities of linear algebras}, Trans. Amer. Math. Soc., {\bf 132} (1968), 553--562.

\bibitem{Kurosh} └.├. ╩єЁю°,
{\it ═хрёёюЎшрЄштэ√х ётюсюфэ√х рыухсЁ√ ш ётюсюфэ√х яЁюшчтхфхэш  рыухсЁ}, ╠рЄхь. ёс.,  {\bf 20} (1947), 239--262.

\bibitem{Shirshov54} └.╚. ╪шЁ°ют,
{\it ╧юфрыухсЁ√ ётюсюфэ√ї ъюььєЄрЄштэ√ї ш ётюсюфэ√ї рэЄшъюььєЄрЄштэ√ї рыухсЁ}, ╠рЄхь. ёс., {\bf 34}(76) (1954), 81--88.

\bibitem{Shirshov53} └.╚. ╪шЁ°ют,
{\it ╧юфрыухсЁ√ ётюсюфэ√ї ышхт√ї рыухсЁ}, ╠рЄхь. ёс., {\bf 33}(75) (1953), 441--452.

\bibitem{Witt} E. Witt,
{\it Die Unterringe der freien Lieschen Ringe},
Math. Z., {\bf 64} (1956), 195--216.

\bibitem{Mikhalev} A.A. ╠шїрыeт,
{\it ╧юфрыухсЁ√ ётюсюфэ√ї ЎтхЄэ√ї ёєяхЁрыухсЁ ╦ш}, ╠рЄхь. чрьхЄъш, {\bf 37}:5 (1985),  653--661.

\bibitem{Shtern} └.╤. ╪ЄхЁэ,
{\it ╤тюсюфэ√х ёєяхЁрыухсЁ√ ╦ш}, ╤шс. ьрЄ. цєЁэ.,  {\bf 27}:1  (1986),  170--174.

\bibitem{Kukin89} ├.┬. ╩Ё цютёъшї, ├.╧. ╩єъшэ,
{\it ╬ яюфъюы№Ўрї ётюсюфэ√ї ъюыхЎ}, ╤шс. ьрЄ. цєЁэ., {\bf 30}:6 (1989), 87--97.

\bibitem{Kukin91} ├.┬. ╩Ё цютёъшї, ├.╧. ╩єъшэ,
{\it └ыуюЁшЄьшўхёъшх ётющёЄтр ётюсюфэ√ї ъюыхЎ}, ╤шс. ьрЄ. цєЁэ., {\bf 32}:6 (1991), 87--99.

\bibitem{AU17} └.└. └ышьсрхт, ╙.╙. ╙ьшЁсрхт, {\it └тЄюьюЁЇшчь ═рурЄ√ ётюсюфэ√ї эхрёёюЎшрЄштэ√ї рыухсЁ Ёрэур фтр эрф хтъышфют√ьш ъюы№Ўрьш}, ╤шс. ¤ыхъЄЁюэ. ьрЄхь. шчт., 14 (2017), 1279--1288

\bibitem{Nagata} M. Nagata,
{\it On the automorphism group of k[x, y]}, Kinokuniya, Tokyo: Kyoto Univ.(Lect. in Math.), 1972.

\bibitem{Cohn72} P.M. Cohn,
{\it Free rings and their relation}, 2nd Ed.-London: Academic Press., 1972.

\bibitem{Shevelin2} ╠.└. ╪хтхышэ, {\it ═хрсхыхт√ 1-ъюуюьюыюушш ш ёюяЁ цхээюёЄ№ ъюэхўэ√ї яюфуЁєяя т эхъюЄюЁ√ї Ёрё°шЁхэш ї,} ╤шс. ьрЄхь. цєЁэ., {\bf53}:5 (2012),  1166--1177.

\bibitem{NU18} ╨.╞. ═рєЁ√чсрхт, ╙.╙. ╙ьшЁсрхт,
{\it ╤ЄЁєъЄєЁр рьры№урьшЁютрээюую яЁюшчтхфхэш  т уЁєяях ртЄюьюЁЇшчьют ётюсюфэ√ї рыухсЁ ╦ш Ёрэур 3}, ╥хчшё√ фюъырфют ьхцфєэрЁюфэющ ъюэЇхЁхэЎшш ╠└╦▄╓┼┬╤╩╚┼ ╫╥┼═╚▀, (2018), 158--158.
http://www.math.nsc.ru/conference/malmeet/18/maltsev18.pdf

\bibitem{ZSSSh} ╩.└. ╞хтыръют, └.╠. ╤ышэ№ъю, ╚.╧. ╪хёЄръют, └.╚. ╪шЁ°ют, {\it ╩юы№Ўр, сышчъшх ъ рёёюЎшрЄштэ√ь,} ╠.: ═рєър, 1978.--431 ёЄЁ.

\bibitem{Shirshov84} └.╚. ╪шЁ°ют, {\it ╩юы№Ўр ш рыухсЁ√}, ╠.: ═рєър, 1984.--143 ёЄЁ.

\bibitem{MKS} ┬. ╠руэєё, └. ╩рЁЁрё, ─. ╤юышЄ¤Ё,  {\it ╩юьсшэрЄюЁэр  ЄхюЁш  уЁєяя. ╧хЁ. ё рэуы.} ╠.: ═рєър, 1974.--455 ёЄЁ.

\bibitem{Vinberg} ▌.┴. ┬шэсхЁу,
{\it ╩єЁё рыухсЁ√}, ╠.: ╘ръЄюЁшры ╧Ёхёё, 2001.-- 544 ёЄЁ.

\end{thebibliography}
\end{document}